\documentclass[11pt, reqno]{amsart}
\usepackage{amsmath, amsthm, a4, latexsym, amssymb}

\def\@rmrk#1#2{\refstepcounter
    {#1}\@ifnextchar[{\@yrmrk{#1}{#2}}{\@xrmrk{#1}{#2}}}

\makeatletter\@addtoreset{equation}{section}\makeatother

 \newfont{\bfit}{cmbxti10 scaled 2000}
 \newfont{\biggi}{cmr12 scaled 2000}

 \newcommand{\eps}{\varepsilon}

 \newcommand{\N}{\mathbb{N}}

 \newcommand{\prob}{\mathbb{P}}

 \newcommand{\one}{\1}

 \newcommand{\skrig}{{\mathcal G}}
 \newcommand{\skrih}{{\mathcal H}}

 \newcommand{\skril}{{\mathcal L}}
 \newcommand{\skrim}{{\mathcal M}}
 \newcommand{\skriw}{{\mathcal W}}

 \newcommand{\skrip}{{\mathcal P}}

 \newcommand{\heap}[2]{\genfrac{}{}{0pt}{}{#1}{#2}}
 \newcommand{\sfrac}[2]{\mbox{$\frac{#1}{#2}$}}

\def\1{{\mathchoice {1\mskip-4mu\mathrm l}      
{1\mskip-4mu\mathrm l}
{1\mskip-4.5mu\mathrm l} {1\mskip-5mu\mathrm l}}}

\newcommand{\eq}{\begin{equation}}
\newcommand{\en}{\end{equation}}

{\nopagebreak {\hspace*{\fill}\rule{2mm}{2mm}}\\ }

{\nopagebreak {\hspace*{\fill}\rule{2mm}{2mm}}\\ }

\renewcommand{\subsection}{\secdef \subsct\sbsect}
\newcommand{\subsct}[2][default]{\refstepcounter{subsection}
\vspace{0.15cm}
{\flushleft\bf \arabic{section}.\arabic{subsection}~\bf #1  }
\nopagebreak\nopagebreak}
\newcommand{\sbsect}[1]{\vspace{0.1cm}\noindent
{\bf #1}\vspace{0.1cm}}

\newtheorem{theorem}{Theorem}[section]
\newtheorem{lemma}[theorem]{Lemma}
\newtheorem{cor}[theorem]{Corollary}

\newtheoremstyle{thm}{1.5ex}{1.5ex}{\itshape\rmfamily}{}
{\bfseries\rmfamily}{}{2ex}{}

\newtheoremstyle{rem}{1.3ex}{1.3ex}{\rmfamily}{}
{\itshape\rmfamily}{}{1.5ex}{}
\theoremstyle{rem}
\newtheorem{remark}{{\slshape\sffamily Remark}}[]

\refstepcounter{subsection}

\def\thebibliography#1{\section*{References}
  \list%
  {\arabic{enumi}.}
    {\settowidth\labelwidth{[#1]}\leftmargin\labelwidth
    \advance\leftmargin\labelsep
    \parsep0pt\itemsep0pt
    \usecounter{enumi}}
    \def\newblock{\hskip .11em plus .33em minus .07em}
    \sloppy                   
    \sfcode`\.=1000\relax}

\begin{document}
\title[Some large deviation  results  for near intermediate random geometric  graphs]
{\Large  Some large deviation  results for near intermediate  random
geometric graphs }

\author[Kwabena Doku-Amponsah ]{}

\maketitle
\thispagestyle{empty}
\vspace{-0.5cm}

\centerline{\sc By Kwabena Doku-Amponsah }
\renewcommand{\thefootnote}{}
\footnote{\textit{AMS Subject Classification:} 60F10, 05C80}
\footnote{\textit{Keywords: }  Random geometric graph, Erd\H
os-R\'enyi graph, coloured random  geometric graph, typed graph,
joint large deviation principle, empirical pair measure, empirical
measure, degree distribution, entropy, relative entropy, isolated
vertices .}
\renewcommand{\thefootnote}{1}

\centerline{\textit{University of Ghana}}

\begin{quote}{\small }{\bf Abstract.}
 We find large
deviation principles for the \emph{degree distribution and the
proportion of isolated vertices } for the \emph{near intermediate}
random  geometric graph models on $n$ vertices placed uniformly in
$[0,1]^d,$  for  $d\in\N.$ In the course of the proof of  these
large deviation results we find joint large deviation principle for
the \emph{ empirical locality  measure } of the coloured random
geometric graphs,(Canning \& Penman, 2003).
\end{quote}\vspace{0.3cm}

\vspace{0.3cm}
\section{Introduction}
In this article we study random  graph model, the  random geometric
graph RGG,  where
 $n$  vertices or  nodes or points are placed uniformly at random in $[0, 1]^d,$
and   any two points distance at most $r_n$ apart are connected. See
(Penrose,  2003). The connectivity radius $r_n$ plays similar role
as the connection probability $p_n$ in  the Erd\H{o}s-R\'{e}nyi
graph  model. Several large  deviation  results  about  the
Erd\H{o}s-R\'{e}nyi graph have  been  established  recently.See
(O'Connell ,1998), (Biggins and Penman, 2009), (Doku-Amponsah and
Moerters, 2010), ( Doku-Amponsah, 2006), (Bordenave and Caputo,
2013), (Mukherjee, 2013)  and (Doku-Amponsah, 2014[a]).

Until  recently few  or  no  large  deviation  result  about  the
degree  distribution  of  the RGG have  been  found. Doku-Amponsah
(2014[b]) proved  some large  deviation principle  for  the degree
distribution of the classical  Erd\H{o}s-R\'{e}nyi graph, where $n$
points are uniformly chosen  in $[0,1]^d$  and  $\lambda_n$ edges
are randomly inserted  among  the  points.

This  article presents a  full large deviation  principle (LDP) for
the empirical degree measure and the proportion of  isolated
vertices of near intermediate RGG presented. Specifically, we  prove
an  LDP for the degree measure  of  the coloured RGG. We Refer to
(Doku-Amponsah and Moerters) for  similar result  for  the
Erd\"{o}-Renyi graphs. From the LDP for  the  empirical degree
measure ; we derive an LDP for the proportion of isolated vertices.
See, O'Connell \cite{OC98} for similar result for the Erd\"{o}-Renyi
graphs.

In  the  course  of  the  proofs  of  this  LDP  we  obtain joint
the empirical pair measure and the \emph{empirical locality
measures} for coloured RGGs. Refer to ( Doku-Amponsah \&
 Moerters, 2010) or (Doku-Amponsah, 2006) for similar results for the coloured
random graphs.

We note  that physical quantities such as the degree distribution,
number of edges per vertex and the proportion of isolated vertices
of RGGs are crucial for understanding many biological systems.

In  the remainder of the paper we state and prove  our  LDP results.
 In Section~\ref{mainresults}  we   state our LDPs,
 Theorem~\ref{Rgg}, Corollary~\ref{RRiv}and
 Theorem~\ref{main2}.  In Section~\ref{proofmain2} we combine (Doku-Amponsah,Theorem~2,1, 2014[b])
 and (Doku-Amponsah,Theorem~2.1, 2014[c]) to obtain the
Theorem~\ref{main2}, using the setup and result of (Biggins, 2004)
to `mix' the LDPs. The paper concludes with the proofs of our  main
results Theorem~\ref{ERdd} and Corollary~\ref{RRiv} which are given
in Section~\ref{proofcorollaries}.

\section{Statement of the results}\label{mainresults}
\subsection{Large  deviations  results for  the random  geometric
graphs.}\label{secsub1} The  RGG is obtained when we sample points
$W_1,...,W_n$  at independently according to the uniform probability
distribution on $[0,\,1]^d,$ for $d\ge 2$ and  given a fixed $r_n>0$
we connect $W_i,W_j$ $(i\not=j)$ if
$$\|W_i-W_j\|\le r_n.$$ See \cite{Pen03}. Various cases of the graph can be
  describe in  terms  of  the  quantity $nr^d,$  which is   a
  measure of  the  average  degree of  the  graph.  See, MCDiarmid  and M\"uller
  \cite{MM05}.  Our  main  aim  in  this  article  is  to present  and
   prove  LDPs   for the empirical degree measure and  the  proportion of isolated nodes to  the
  number  of  vertices  of the RGG  when  the  connectivity
  radius satisfies  $nr^d\to \alpha,$  for  $\alpha>0.$  Thus, we  consider
  the  near  intermediate  case.

The first
  theorem  in  this  subsection is the  LDP  for  the  degree
  distribution  of  the  RGGs.  We  assume  $d\ge 2$ is finite and
  write
  $$\rho(d)=\sfrac{\pi^{d/2}}{\Gamma\big(\sfrac{(d+2)}{2}\big)},$$
  where  $Gamma$  is  the  gamma  function.

\begin{theorem}\label{ERdd}\label{Rgg}
Suppose $D$ is the degree distribution of  the random graph
$\skrig(n,r_n),$  where the connectivity radius $r_n\in (0,1]$
satisfies $n r_n^d \to c \in (0,\infty)$. Then ,as $n\to\infty$, $D$
satisfies an LDP in the space $\skrip(\N \cup \{0\})$ with good rate
function
\begin{equation}\label{randomg.ratedeg}
\begin{aligned}
\eta_1(\delta)= \left\{ \begin{array}{ll}  \frac 12 \, \langle
\delta \rangle\, \log \Big( \sfrac{\langle \delta\rangle}{\rho(d)c}
\Big)
 - \frac 12 \, \langle \delta \rangle  +  \sfrac {\rho(d)c}{2} +H (d\,\|\,q_{\langle \delta\rangle}),
 & \mbox { if $\langle \delta\rangle< \infty,$ }\\[2mm]
\infty & \mbox { if $\langle \delta\rangle= \infty,$ }
\end{array}\right.
\end{aligned}
\end{equation}
where $q_{x}$ is a Poisson distribution with parameter~$x$ and
$\langle \delta\rangle:= \sum_{m=0}^{\infty}m\delta(m)$.
\end{theorem}

Next  we give  a  similar result as in O'Connell~\cite{OC98}, the
LDP for the  proportion of isolated vertices of the  RGG.

 \begin{cor}\label{ERiv}\label{RRiv}
Suppose $D$ is the degree distribution of  the random graph
$\skrig(n,r_n),$ where the connectivity radius $r_n\in (0,1]$
satisfies $n r_n^d \to c \in (0,\infty)$. Then, as $n\to\infty$, the
proportion of isolated vertices, $D(0)$ satisfies an LDP in $[0,1]$
with good rate function $$\xi_1(y)=y \log y + \rho(d) c y(1-y/2) -
(1-y) \big[ \log\big (\sfrac {\rho(d)c}{a}\big) -
 \sfrac{(a- \rho(d) c (1-y))^2}{2\rho(d)c(1-y)} \big] \, ,$$
 where $a=a(y)$ is the unique positive solution of $1-e^{-a}=\frac {\rho(d)c}{a}\, (1-y)$.
 \end{cor}

 From  Lemma~\ref{RRiv} we  deduce  that on a typical random  geometric graphs  the
 number  of isolated  vertices will  grow  like  $n
 e^{-\rho(d)c}.$ Thus,  as  $n\to\infty,$  the  number  of isolated
 vertices in the R.G graphs  converges  to $n
 e^{-\rho(d)c}$  in probability.
 In our  last  theorem  in  this  subsection  we  give  the
 LDP for  the  proportion  of  edges to the  number
  of  vertices  of  the  R.G.

\subsection{Large deviation principles for empirical  measures of  the  coloured
random geometric graphs.}  In this subsection  we shall look at  a
more general model of random  geometric graphs, the coloured RGGs in
which the connectivity radius depends on  the type or colour or
symbol or  spin of the  nodes. The empirical pair measure and the
empirical locality  measure are our main object of study.

Given   a probability measure $\nu$ on $\skriw$ and a function
$r_n\colon\skriw\times\skriw\rightarrow (0,1]$ we may define the
{\em randomly coloured random geometric graph} or simply
\emph{coloured random geometric graph}~$X$ with $n$ vertices as
follows: Pick vertices  $W_1,...,W_n$  at  random  independently
according to the uniform distribution on $[0,\,1]^2.$ Assign to each
vertex $W_j$ colour $X(W_j)$ independently according to the {\em
colour law} $\mu.$ Given the colours, we join any two vertices
$W_i,W_j$,$(i\not=j)$ by an edge independently of everything else,
if $$\|W_i-W_j\|\le r_n\big[X(W_i),X(W_j)\big].$$ In this  article
we  shall  refer to $r_n(a,b),$   for $a,b\in\skriw$ as a connection
radius,  and always consider
$$X=((X(W_i),X(W_j))\,:\,i,j=1,2,3,...,n),E)$$ under the joint law of
graph and colour. We interpret $X$ as  coloured RGG with vertices
$Y_1,...,Y_n$ chosen at  random uniformly   and independently from
the vertices space $[0,1]^2.$ For  the  purposes of  this  study we
restrict ourselves to  the near intermediate cases .i.e. the
connection radius $r_n$ satisfies the condition $n r_n^d(a,b) \to
C(a,b)$ for all $a,b\in \skriw$, where $C\colon\skriw^2\rightarrow
[0,\infty)$ is a symmetric function, which is not identically equal
to zero.

For any finite or countable set $\skriw$ we denote by
$\skrip(\skriw)$ the space of probability measures, and by
$\tilde\skrip(\skriw)$ the space of finite measures on $\skriw$,
both endowed with the weak topology.  By  convention we write
$\N=\{0,1,2,...\}.$

We associate with any coloured graph $X$ a probability measure, the
\emph{empirical colour measure}~$\skril^1\in\skrip(\skriw)$,~by
$$\skril_X^{1}(a):=\frac{1}{n}\sum_{j=1}^{n}\delta_{X(W_j)}(a),\quad\mbox{ for $a_1\in\skriw$, }$$
and a symmetric finite measure, the \emph{empirical pair measure}
$\skril_X^{2}\in\tilde\skrip_*(\skriw^2),$ by
$$\skril_X^{2}(a,b):=\frac{1}{n}\sum_{(i,j)\in E}[\delta_{(X(W_i),X(W_j))}+
\delta_{((X(W_j),X(W_i))}](a,b),\quad\mbox{ for $(a, b)\in\skriw^2$.
}$$ The total mass $\|\skril_X^2\|$ of the empirical pair measure is
$2|E|/n$.  Finally we define a further probability measure, the
\emph{empirical neighbourhood measure}
$\skrim_X\in\skrip(\skriw\times\N)$, by
$$\skrim_X(a,\ell):=\frac{1}{n}\sum_{j=1}^{n}\delta_{(X(W_j),L(W_j))}(a,\ell),\quad
\mbox{ for $(a,\ell)\in\skriw\times\N$, }$$ where
$L(v)=(l^{v}(b),\,b\in\skriw)$ and $l^{v}(b)$ is the number of
vertices  of colour $b$  connected to vertex $v$.

For any $\mu\in\skrip(\skriw\times\N^{\skriw})$we  denote by $\mu_1$
the $\skriw-$ marginal of $\mu$ and for  every
$(b,a)\in\skriw\times\skriw,$ let $\mu_2$  be the law  of  the pair
$(a,l(b))$ under  the  measure $\mu.$ Define the  measure (finite),
$\langle\mu(\cdot,\ell),\,l(\cdot)\rangle\in\tilde\skrip(\skriw\times\skriw)$
by
$$\skrih_2(\mu)(b,a):=
\sum_{l(b)\in\N}\mu_2(a,l(b))l(b), \quad\mbox{ for $a,b\in\skriw$}$$
and  write $\skrih_1(\mu)=\mu_1.$ We define the function  $\skrih
\colon \skrip(\skriw\times\N^{\skriw}) \to \skrip(\skriw) \times
\tilde\skrip(\skriw \times \skriw)$ by
$\skrih(\mu)=(\skrih_1(\mu),\skrih_2(\mu))$ and note that
$\skrih(\skrim_X)=(\skril_X^1, \skril_X^2).$ Observe  that
$\skrih_1$ is a continuous function  but $\skrih_2$ is
\emph{discontinuous} in the weak topology. In particular, in  the
summation  $\displaystyle \sum_{l(b)\in\N}\mu_2(a,l(b))l(b)$ the
function $l(b)$ may be unbounded and  so  the  functional
$\displaystyle \mu\to\skrih_2(\mu)$ would not be continuous in the
weak topology. We call a pair of measures
$(\varpi,\mu)\in\tilde{\skrip}(\skriw\times\skriw)\times\skrip(\skriw\times\N^{\skriw})$
\emph{sub-consistent} if
\begin{equation}\label{consistent}
\skrih_2(\mu)(b,a) \le \varpi(b,a), \quad\mbox{ for all
$a,b\in\skriw,$}
\end{equation}
and \emph{consistent} if equality holds in \eqref{consistent}. For a
measure $\varpi\in\tilde\skrip_*(\skriw^2)$ and a measure
$\omega\in\skrip(\skriw)$, define
$${{\mathfrak H}_C^{d}}(\varpi\, \| \, \omega ):=
H\big(\varpi\,\|\,\rho(d)C\omega\otimes\omega\big)+\rho(d)\|
C\omega\otimes\omega \| -\|\varpi\|\, ,$$ where the measure
$C\omega\otimes\omega\in\tilde\skrip(\skriw\times\skriw)$ is defined
by $C\omega\otimes\omega(a,b)=C(a,b)\omega(a)\omega(b)$ for
$a,b\in\skriw$. It is not hard to see that $\mathfrak
H_C^{d}(\varpi\,\|\,\omega)\ge 0$ and equality holds if and only if
$\varpi=\rho(d) C\omega\otimes\omega$. For every
$(\varpi,\mu)\in\tilde\skrip_*(\skriw\times\skriw) \times
\skrip(\skriw\times\N)$ define a probability measure
$Q=Q[\varpi,\mu]$ on $\skriw\times\N$ by
$$Q(a\,,\,\ell):=\mu_{1}(a)\prod_{b\in\skriw} e^{-\frac{\varpi(a,b)}{\mu_1(a)}} \,
\frac{1}{\ell(b)!}\,
\Big(\frac{\varpi(a,b)}{\mu_1(a)}\Big)^{\ell(b)}, \quad\mbox{for
$a\in\skriw$, $\ell\in\N$} .$$ We now  state the principal theorem
in  this section the LDP for the empirical pair measure and the
empirical neighbourhood measure.

\begin{theorem}\label{randomg.LDM}\label{main2}
Suppose that $X$ is a coloured RGG  graph  with colour law $\mu$ and
connection radii $ r_n\colon\skriw\times\skriw\rightarrow[0,1]$
satisfying $n r_n^d(a,b) \to C(a,b)$ for some symmetric function
$C\colon\skriw\times\skriw\rightarrow [0,\infty)$ not identical to
zero. Then, as $n\rightarrow\infty,$ the pair
$(\skril_X^2,\,\skrim_X)$ satisfies an LDP in
$\tilde{\skrip}_*(\skriw\times\skriw)\times\skrip(\skriw\times\N)$
with good rate function
$$\begin{aligned}
J(\varpi,\mu)=\left\{
\begin{array}{ll}H(\mu\,\|\,Q)+H(\mu_1\,\|\,\nu)+
\sfrac{1}{2}\, {\mathfrak H}_C^{d}(\varpi\,\|\,\mu_1 ) & \mbox {if $(\varpi,\mu)$ consistent  and  $\mu_1=\varpi_2,$   }\\
\infty & \mbox{otherwise.}
\end{array}\right.
\end{aligned}$$
\end{theorem}

\begin{remark}

Note that on  typical coloured RGG graph
 we  have, $\omega=\mu_1,$
$\varpi= \rho(d) C \,\mu \otimes \mu$ and
$$\mu(a,\ell)=\nu(a)\prod_{b\in\skriw}e^{-\rho(d) C(a,b)\nu(b)}\, \frac{(\rho(d)C(a,b)\nu(b))^{\ell(b)}}{\ell(b)!},\qquad\mbox{for
all $(a,\ell)\in\skriw\times\N$} .$$ This is the law of a pair
$(a,\ell)$ where $a$ is distributed according to $\mu$ and, given
the value of $a$, the random variables $\ell(b)$ are independently
Poisson distributed with parameter $\rho(d) C(a,b)\nu(b)$.  Hence,
as  $n\to\infty,$  the  empirical neigbourhood  measure
$\skrim_X(a,\ell)$ converges to  $\mu(a,\ell)$  in probability.
\end{remark}

\section{Proof of Theorem~\ref{main2}}\label{proofmain2}
For any $n\in\N$ we define
$$\begin{aligned}
\skrip_n(\skriw) & := \big\{ \omega\in \skrip(\skriw) \, : \, n\omega(a) \in \N \mbox{ for all } a\in\skriw\big\},\\
\tilde \skrip_n(\skriw \times \skriw) & := \big\{ \varpi\in
\tilde\skrip_*(\skriw\times\skriw) \, : \, \sfrac
n{1+\one\{a=b\}}\,\varpi(a,b) \in \N  \mbox{ for all } a,b\in\skriw
\big\}\, .
\end{aligned}$$

We denote by
$\Theta_n:=\skrip_n(\skriw)\times\tilde{\skrip}_n(\skriw\times\skriw)$
and
$\Theta:=\skrip(\skriw)\times\tilde{\skrip}_*(\skriw\times\skriw)$.
With
$$\begin{aligned}
P_{(\omega_n, \varpi_n)}^{(n)}(\mu_n) & := \prob\big\{\skrim_X=\mu_n \, \big| \, \skrih(\skrim_X)=(\omega_n,\varpi_n)\big\}\, ,\\
P^{(n)}(\omega_n,\varpi_n) & :=
\prob\big\{(\skril_X^1,\skril_X^2)=(\omega_n,\varpi_n)\big\}
\end{aligned}$$

the joint distribution of $\skril_X^1, \skril_X^2$ and $\skrim_X$ is
the mixture of $P_{(\omega_n, \varpi_n)}^{(n)}$ with
$P^{(n)}(\omega_n,\varpi_n)$ defined as
\begin{equation}\label{randomg.mixture}
d\tilde{P}^n(\omega_n, \varpi_n, \mu_n):= dP_{(\omega_n,
\varpi_n)}^{(n)}(\mu_n)\, dP^{(n)}(\omega_n, \varpi_n).\,
\end{equation}

(Biggins, Theorem 5(b), 2004) gives criteria for the validity of
large deviation principles for the mixtures and for the goodness of
the rate function if individual large deviation principles are
known. The following three lemmas ensure validity of these
conditions.

\begin{lemma}[Doku-Amponsah, 2014b] \label{randomg.uniexpotightness} The  family of
measures $({P}^n \colon n\in\N)$  is  exponentially tight on
$\Theta$
\end{lemma}

\begin{lemma}[Doku-Amponsah \& Moerters, 2010] \label{randomg.uniexpotightness} The  family of
measures $(\tilde{P}^n \colon n\in\N)$  is  exponentially tight on
$\Theta\times\skrip(\skriw\times\N).$
\end{lemma}

Define the function
$$\tilde{J}\colon{\Theta}\times\skrip(\skriw\times\N)\rightarrow[0,\infty],
\qquad\tilde{J}((\omega,\varpi),\,\mu)=\tilde{J}_{(\omega,\varpi)}(\mu),$$
where

\begin{align}\label{randomg.rateLDprob}
\tilde{J}_{((\varpi,\omega)}(\mu)=\left\{
\begin{array}{ll}H(\mu\,\|\,Q_{poi}) & \mbox
  {if  $(\omega,\mu)$  is  consistent and $\mu_1=\omega_2$ }\\
\infty & \mbox{otherwise.}
\end{array}\right.
\end{align}

\begin{lemma}[Doku-Amponsah \& Moerters, 2010]\label{randomg.convexgoodrate}
$\tilde{J}$ is lower semi-continuous.
\end{lemma}

By (Biggins, Theorem~5(b), 2004) the two previous lemmas and the
large deviation principles we have established in (Doku-Amponsah,
Theorem~2.1, 2014b) and (Doku-Amponsah, Theorem~2.1, 2014c) ensure
that under $(\tilde{P}^n)$ the random variables $(\omega_n,
\varpi_n, \mu_n)$ satisfy a large deviation principle on
$\skrip(\skriw) \times \tilde\skrip_*(\skriw\times\skriw)
\times\skrip(\skriw\times\N)$ with good rate function
$$\hat{J}(\omega, \varpi, \mu)= \left\{ \begin{array}{ll}
H(\omega\,\|\,\nu) + \sfrac12\,{\mathfrak H}_C^d(\varpi\,\|\,\omega)
+
H(\mu\,\| Q_{poi})\, , & \mbox{ if $(\varpi,\mu)$  is  consistent and $\mu_1=\varpi_2,$} \\
 \infty\, , & \mbox{ otherwise.} \\ \end{array}\right. $$
By projection onto the last two components we obtain the large
deviation principle as stated in Theorem~\ref{randomg.LDM} from the
contraction principle, see e.g. (Dembo et al.,1998, Theorem~4.2.1).

\section{Proof of Theorem~\ref{ERdd}~and~Corollary~\ref{RRiv}}\label{proofcorollaries}

We derive the theorems from Theorem~\ref{main2} by applying the
contraction principle, see e.g.~(Dembo \& Zeitouni,Theorem~4.2.1,
1998). In
 fact Theorem~\ref{main2} and the contraction principle imply a large
deviation principle for $D$. It just remains to simplify the rate
functions.

\subsection{Proof of Theorem~\ref{ERdd}.} Note  that,  in the case of an uncoloured
RGG graphs, the function $C$ degenerates to a constant~$c$,
$L^2=|E|/n\in[0,\infty)$ and $M=D\in\skrip(\N\cup\{0\})$.
Theorem~\ref{main2}  and the contraction principle imply a large
deviation principle for $D$ with good rate function
$$\begin{aligned} \eta_1(\delta) & =\inf\big\{J(x,\delta) \colon x\ge 0 \big\}
=\inf\big\{ H(\delta\,\|\,q_{x})+\sfrac1 2 x\log x -\sfrac12 x\log
\rho(d)c+ \sfrac{1}{2}\,\rho(d)c-\sfrac12 x \colon \langle
\delta\rangle \le x\,\big\},
\end{aligned}$$ which is to be understood as infinity if $\langle
d\rangle$ is infinite. We denote by $\eta^{x}(\delta)$ the
expression inside the infimum. For any $\eps>0$, we have
$$\begin{aligned}
\eta^{\langle \delta\rangle+\eps}(\delta)-\eta^{\langle
\delta\rangle}(\delta) & =\sfrac{\eps}{2} +\sfrac{\langle
\delta\rangle-\eps}{2}\log\sfrac{\langle \delta\rangle}{\langle
\delta\rangle+\eps}+\sfrac{\eps}{2}\log\sfrac{\langle
\delta\rangle}{\rho(d)c} \ge \sfrac{\eps}{2} +\sfrac{\langle
\delta\rangle-\eps}{2}\, \big(\sfrac{-\eps}{\langle
\delta\rangle}\big)+\sfrac{\eps}{2}\log\sfrac{\langle
\delta\rangle}{\rho(d)c}
>0,
\end{aligned}$$
so that the minimum is attained at $x=\langle\delta\rangle$.

 \subsection{Proof of Corollary~\ref{RRiv}.} Corollary~\ref{ERiv} follows from Theorem~\ref{ERdd} and the
contraction principle applied to the continuous linear map
$G\colon\skrip(\N\cup\{0\})\rightarrow[0,\,1]$ defined by
$G(\delta)=\delta(0).$ Thus, Theorem~\ref{ERdd} implies the large
deviation
 principle for $G(D)=W$ with the good rate function
$\xi_1(y)=\inf\{\eta_1(\delta) \colon \delta(0)=y, \langle
\delta\rangle < \infty\}.$ We  recall  the  definition  of
$\eta^{x}$ and  observe that $\xi_2(y)$ can be expressed as
$$\xi_1(y)=\inf_{b\ge
0}\inf_{\heap{d\in\skrip(\N\cup\{0\})}{\delta(0)=y,\,
\rho(d)c\langle \delta\rangle=b^2}} \Big\{\sfrac{1}{2}c+y\log
y+\sfrac{b^2}{2\rho(d)c} +\sum_{k=1}^{\infty}\delta(k)\log
\sfrac{\delta(k)}{q_{b}(k)}-b(1-y)\Big\}.$$ Now, using Jensen's
inequality, we have that
\begin{equation}\label{randomg.Jensen}\sum_{k=1}^{\infty}\delta(k)\log
\sfrac{\delta(k)}{q_{b}(k)}\ge
(1-y)\log\sfrac{(1-y)}{(1-e^{-b})},\end{equation}
 with equality if $\delta(k)=\sfrac{(1-y)}{(1-e^{-b})}q_{b}(k),$ for all
$k\in\N.$  Therefore, we have the inequality
$$\inf\big\{\eta(\delta) \colon \delta(0)=y,
\langle \delta\rangle <\infty\big\} \ge \inf\big\{
\sfrac{1}{2}c+y\log y+\sfrac{b^2}{2\rho(d)c}
+(1-y)\log\sfrac{(1-y)}{(1-e^{-b})} -b(1-y) \colon b\ge 0 \big\}.$$
Let $y\in[0,\,1].$  Then, the equation $a(1-e^{-a})=\rho(d)c(1-y)$
has a unique positive solution. Elementary calculus shows that the
global minimum of
 $b\mapsto\sfrac{1}{2}\rho(d)c+y\log y+\sfrac{b^2}{2\rho(d)c} +(1-y)\log\sfrac{(1-y)}{(1-e^{-b})}
-b(1-y)$ on $(0,\infty)$ is attained at the value $b=a$,
 where $a$ is the positive solution of our equation.
We obtain the form of $\xi$  in Corollary~\ref{ERiv} by observing
that
$$\sfrac{a(y)^2+(\rho(d)c)^2 -2\rho(d)ca(y)\big(1-y\big)}{2\rho(d)c}=\sfrac{\rho(d)cy}{2}\big(2-y\big)+\sfrac{1}{2\rho(d)c}\big(a(y)-\rho(d)c(1-y)\big)^2.$$

\bigskip


\end{document}